\documentclass[final,11pt]{article}
\usepackage{amsmath}
\usepackage{amsfonts}
\usepackage{amsthm}

\def\R{\mathbb{R}}

\begin{document}

\title{Solution to the Pompeiu problem and the related symmetry problem}

\author{A. G. Ramm\\
\small Mathematics Department\\
\small Kansas State University, Manhattan, KS 66506, USA\\
\small email: ramm@math.ksu.edu}
\date{}
\maketitle

\begin{abstract}

Assume that $D\subset \R^3$ is a bounded domain
with $C^1-$smooth boundary. Our result is:

{\bf Theorem 1.} {\em If $D$  has $P-$property, then $D$ is a ball.}

Four equivalent formulations of the Pompeiu problem are discussed.

A domain $D$ has $P-$property if there exists an $f\neq 0$, $f\in L^1_{loc}(\mathbb{R}^3)$
such that $\int_{D}f(gx+y)dx=0$ for all $y\in \mathbb{R}^3$ and all $g\in SO(2)$,
where $ SO(2)$ is the rotation group.

The result obtained concerning the related symmetry problem is:

{\bf Theorem 2.} {\em If $(\nabla^2 +k^2)u=0$ in $D$, $u|_S=1$, $u_N|_S=0$, and $k>0$ is a constant, then $D$
is a ball.}

\end{abstract}



\section{Introduction} \label{Introduction}

The modern formulation of the Pompeiu problem can be given in
several equivalent ways. The history of this problem goes
back to 1929, see \cite{P}.
 Brief historical remarks are at the end of our paper.
Some known results and references
about this problem can be found in \cite{R470}, Chapter 11,
\cite{R363}, \cite{R629},  \cite{W}, \cite{Z}. Our bibliography is incomplete.
  Our paper is essentially self-contained.

We assume throughout that $D\subset \R^3$ is a bounded domain
 and its boundary $S$ is $C^1-$smooth.
 It is known
that under the assumptions of any of the {\bf Formulations 1-4}, stated below,
the boundary $S$ is real-analytic, see \cite{W}.

By $\chi$ the characteristic
function of $D$ and by $\tilde{\chi}(\xi):=\int_De^{i\xi\cdot x}dx$
its Fourier transform are denoted, respectively. The vector
$\xi=k\alpha$, where $k>0$ is the length of $\xi$, $\alpha\in S^{2}$
is a unit vector in $\R^3$, $S^{2}$ is a unit sphere in $\R^3$, and
$\xi\cdot x=(\xi, x)$ is the dot product in $\R^3$, $[\xi, x]$ is the cross product.

One of the modern formulations of the Pompeiu problem is the
following (\cite{R470}, \cite{R629}):

{\bf Formulation 1.} {\it Assume that $f\in L^1_{loc}$, $f\not\equiv 0$, and
\begin{equation}\label{eq:1}
\int_Df(y+gx)dx=0, \quad \forall y\in \R^3, \quad \forall g,
\end{equation}
where $g$ is an arbitrary rotation. Prove that $D$ is a
ball.}

{\bf Formulation 2.} {\it Prove that if
\begin{equation}\label{eq:3}
\tilde{\chi}(k\alpha)=0
\end{equation}
for all $\alpha\in S^2$ and a fixed number $k>0$, then $D$ is a ball.}

Let us give other equivalent formulations of the Pompeiu problem.
Denote by $N$ the unit normal to $S$ pointing out of $D$.

{\bf Formulation 3.} {\it Suppose that $k>0$ is fixed and the following problem
\begin{equation}\label{eq:2}
(\nabla^2+k^2)u=c \quad in \quad D, \quad u|_S=u_N|_S=0,\quad c=const\neq 0,
\end{equation}
has a solution. Prove that $D$ is a ball.}

Without loss of generality one may assume that the constant $c$ in (\ref{eq:2})
equals to $1$. We assume this below.

Clearly, {\bf Formulation 3} is equivalent to the following symmetry problem, related to the Pompeiu problem:

{\bf Formulation 4.} {\it Suppose that $k>0$ is fixed and the following problem
\begin{equation}\label{eq:2'}
(\nabla^2+k^2)u=0 \quad in \quad D, \quad u|_S=const\neq 0,\quad   u_N|_S=0,
\end{equation}
has a solution. Prove that then $D$ is a ball.}

{\bf Formulation 4} gives a symmetry problem that has been unsolved for decades.

If, as we assumed above, the constant $c=1$ in {\bf Formulation 3} then
the $const$ in   {\bf Formulation 4} equals to $-k^{-2}$.

{\bf Formulation 3} with $c=1$ and {\bf Formulation 4} are equivalent. To prove this,  just make in
(\ref{eq:2}) the substitution $u=v+k^{-2}$.

The equivalence of {\bf Formulations 1, 2, 3 and 4} is  proved below,
 see also \cite{R470}, Chapter 11, and \cite{R629}. In our
proof ${\bf 1}$ stands for {\bf Formulation 1}, etc.

\hspace{3mm}

{\bf Proof of the equivalence of Formulations 1, 2, 3, 4:}

\hspace{2mm}

${\bf 2\Rightarrow 3.}$ If an entire function of exponential type
$\tilde{\chi}(\xi)$ vanishes on the irreducible algebraic variety
$\xi^2=k^2$, then the function
$\tilde{u}:=\tilde{\chi}(\xi)(\xi^2-k^2)^{-1}$ is also entire and of the
same exponential type. Its Fourier transform $u(x)$ solves problem
(\ref{eq:2}). The function $u$ is defined in all of $\mathbb{R}^3$, $u\in H^2_{loc}(\mathbb{R}^3)$ and
is compactly supported by the Paley-Wiener theorem. Therefore, by the unique continuation theorem,
 $u=0$ in $D':=\mathbb{R}^3\setminus D$. This and the embedding theorem imply the boundary
 conditions in {\bf Formulation 3}.  \hfill $\Box$

\hspace{2mm}

${\bf 3\Rightarrow 2.}$ If (\ref{eq:2}) holds, then multiply (\ref{eq:2})
by $e^{ik\alpha \cdot x}$, $\alpha\in S^2$ is arbitrary, integrate over $D$, and
then by parts, using the boundary conditions (\ref{eq:2}) and the equation
$(\nabla^2+k^2)e^{ik\alpha \cdot x}=0$. This yields (\ref{eq:3}). \hfill $\Box$

\hspace{2mm}

${\bf 1\Rightarrow 2.}$ Take the Fourier transform (in the distributional
sense) of (\ref{eq:1}) and get
$\overline{\tilde{\chi}(g^{-1}\xi)}\tilde{f}(\xi)=0$ for all $g$, where the over-line
stands for  complex conjugate. Therefore,
$supp \tilde{f}=\cup_{k}C_k$, where
$$C_k:=\{\xi: \xi^2=k^2, \tilde{\chi}|_{\xi^2=k^2}=0\}, \quad k>0,$$
 and the set $\{k\}$ is a discrete
set of positive numbers since $\tilde{\chi}$ is an entire function.
Thus, there is a $k>0$ such that (\ref{eq:3}) holds. \hfill $\Box$

\hspace{2mm}

${\bf 3\Rightarrow 1.}$ If (\ref{eq:2}) holds, then one extends $u$ to $\mathbb{R}^3$
by setting $u=0$ in $D'$, and takes Fourier transform of this $u$ assuming $c=1$.
This yields $\tilde{u}(\xi)(k^2-\xi^2)=\tilde{\chi}(\xi)$.
Let
$\tilde{f}\neq 0$ be supported on $C_k:=\{\xi: \xi^2=k^2, \tilde{\chi}|_{\xi^2=k^2}=0\}$, where $k>0$ is a number. Then,
$\overline{\tilde{\chi}(g^{-1}\xi)}\tilde{f}(\xi)=0$ for any rotation $g$. Taking the inverse
Fourier transform of this relation yields (\ref{eq:1}). \hfill $\Box$

\hspace{2mm}

From the above relations it follows that ${\bf 1}\Longleftrightarrow {\bf 2}\Longleftrightarrow {\bf 3}$. \hfill $\Box$

We have already proved that ${\bf 3}\Longleftrightarrow {\bf 4}$. \hfill $\Box$

Thus, {\bf Formulations 1, 2, 3, 4} are equivalent:

${\bf 1}\Longleftrightarrow {\bf 2}\Longleftrightarrow {\bf 3}\Longleftrightarrow {\bf 4}$.   \hfill$\Box$

\hspace{2mm}

Because {\bf Formulation 1} is equivalent to {\bf Formulations 2, 3} and {\bf 4}
the following  Theorems 1, 2 and 4 will be established if Theorem 3 is proved.

\hspace{2mm}
\newpage

{\bf Theorem 1.} {\em If (1) holds, then $D$ is a ball.}

\hspace{2mm}

{\bf Theorem 2.} {\it If (\ref{eq:3}) holds, then $D$ is a ball.}

\hspace{2mm}

{\bf Theorem 3.} {\it If (\ref{eq:2}) holds, then $D$ is a ball.}

\hspace{2mm}

{\bf Theorem 4.} {\it If (\ref{eq:2'}) holds, then $D$ is a ball.}

\hspace{2mm}

In Section 2 we prove Theorem 3 and, therefore, Theorems 1, 2 and 4.

\section{Proof of Theorem 3.}

  If equation (\ref{eq:2}) holds,
then one derives the following equations (see also \cite{R470}, p. 415, \cite{R382}):
\begin{equation}\label{eq:4}
[s,N]=\sum_{j=1}^J c_ju_{jN}(s), \qquad \forall s\in S, \quad c_j=const,
\end{equation}
where $ c_j$ are constant vectors, $[s,N]$ is the cross product,  $u_j$ are solutions to
\begin{equation}\label{eq:5}
(\nabla^2+k^2)u_j=0, \qquad u_j{|_S}=0, \qquad 1\le j\le J,
\end{equation}
and the system $\{u_j\}_{j=1}^J$ is maximal linearly independent set of solutions to  (\ref{eq:5}) normalized in $L^2(D)$.

We prove formula (\ref{eq:4}) at the end of the paper, see Theorem 6.

{\bf Theorem 5.} {\em If  (\ref{eq:4})-(\ref{eq:5}) hold and $S$ is real-analytic, then $c_j=0$,
$1\le j \le J$, so $[s,N]=0$ for all $s\in S$. Therefore $S$ is a sphere centered at the origin.}

{\em Proof.} If the origin $O$  of the coordinate system is moved to a position, described
by vector $a$, then $s\to s+a$, $N\to N$, $u_{jN}\to u_{jN}$ and $[s,N]\to [s,N]+[a,N]$. The coefficients
 $c_j\to c_j(a)$, where $c_j$ and $c_j(a)$ do not depend on $s$ and $c_j$ do not depend on $a$.
Write equation  (\ref{eq:4}) with $s+a$, take into account that $N(s)$ and $u_{jN}(s)$ do not depend on $a$. This
equation is valid for any $a$:
\begin{equation}\label{eq:4a}
[s+a,N(s)]=\sum_{j=1}^Jc_j(a)u_{jN}(s).
\end{equation}
Subtract from the equation  (\ref{eq:4a}) equation  (\ref{eq:4}) and get
\begin{equation}\label{eq:6}
[a,N(s)]=\sum_{j=1}^Jb_ju_{jN}(s), \qquad b_j:=c_j(a)-c_j.
\end{equation}
Let us prove that the vectors $b_j=0$, $1\le j \le J$.
 Note that $c_j(a)$ depends on $a$ but does not depend on $s$, and $c_j$ does  not depend on $a$ or on $s$. Take
an arbitrary unit vector $q$.
 Denote $B_j:=q\cdot b_j$. Let $p=[q,a]$ and $N_p$ be the projection of $N$ onto $p$.
 Scalar multiply (\ref{eq:6}) by $q$ and get:
 \begin{equation}\label{eq:6a}
N_p(s)=\sum_{j=1}^J B_ju_{jN}(s), \qquad B_j:=q\cdot b_j.
\end{equation}
Choose a closed curve $\mathcal{L}\in S$ such that $N_p(s)|_{\mathcal{L}}=0$. Existence of such a curve
is clear geometrically: this curve consists of the points $s$ of $S$ which belong to the straight lines
tangent to $S$ and parallel to $p$.
Assume that one can  choose $J$ points $s_m\in \mathcal{L}$ such that $det (u_{jN}(s_m))\neq 0$.
In Lemma 3 it is proved that this determinant is not equal to $0$ for almost all $s_m\in S$ and is an analytic
function of $s_m$ on $S$. Since $q$ is an arbitrary unit vector, $p$ has an arbitrary direction
orthogonal to $a$. Therefore, one may assume that $det (u_{jN}(s_m))\neq 0$, because a slight change
of $q$ allows one to slightly move $\mathcal{L}$. So, if on every $\mathcal{L}$ the determinant is equal to zero
for all $s_m$, then this determinant is equal to zero on a set of positive measure on $S$. This contradicts
the conclusion of Lemma 3. Thus, it follows that there exist $J$ points $s_m\in \mathcal{L}$ such that
$det (u_{jN}(s_m))\neq 0$. Therefore,
equation   (\ref{eq:6a}) with  $s=s_m$, $1\le m \le J$, implies that all $B_j=0$.
Since $q$ is arbitrary, equation $B_j=0$ implies that $b_j=c_j(a)-c_j=0$.

Since $c_j(a)=c_j$,  equation (\ref{eq:4a}) can be rewritten as
\begin{equation}\label{eq:4b}
[s+a,N(s)]=\sum_{j=1}^Jc_j u_{jN}(s),
\end{equation}
where $c_j$ do not depend on $a$, so the right side of this equation is bounded as $a\to \infty$.
Consequently, if not all $c_j$ are equal to zero one has a contradiction, since the left side
of  equation (\ref{eq:4b}) grows (for some $s$) as $a\to \infty$ while its right side is bounded
as $a\to \infty$. If all $c_j=0$, then equation (\ref{eq:4b})
yields the relation $[s+a,N(s)]=0$ on $S$. This relation implies that $[s,N(s)]=0$ on $S$ and $[a,N(s)]=0$ on $S$.
If $[s,N(s)]=0$ on $S$, then  $S$ is a sphere by Lemma 1, see below.
If $[a,N(s)]=0$ on $S$, then $a=0$ since $a$ does not depend on $s$ and $N(s)$
depends on $s$. If $a=0$ then the sphere $S$ has its center at the origin.

Theorem 5 is proved. \hfill$\Box$

{\bf Lemma 1.} (\cite{R470}, p. 336)  {\em If $S$ is a smooth closed surface in $\mathbb{R}^3$ such that $[s,N]=0$ for all $s\in S$, then $S$ is a sphere.}

{\em Proof.} Let $s=s(p,q)$ be a parametric equation of $S$. Then the normal $N$ is proportional to $[s_p,s_q]$, so
\begin{equation}\label{eq:7}
0=[s,[s_p,s_q]]=s_p s\cdot s_q-s_q s\cdot s_p.
\end{equation}
Vectors $s_p$ and $s_q$ are linearly independent at regular points of $S$. Since $S$ is analytic all its points are regular.
Therefore, equation  (\ref{eq:7}) implies that $s\cdot s_q=0$ and $s\cdot s_p=0$. Consequently, $s\cdot s=const$. This means that
$S$ is a sphere. Lemma 1 is proved. \hfill$\Box$

{\bf Lemma 2.} {\em If  equations (\ref{eq:5}) hold and the system $\{u_j(x)\}_{j=1}^J$ is linearly independent in
$L^2(D)$, then the system $\{u_{jN}(s)\}_{j=1}^J$ is linearly independent in $L^2(S)$.}

{\em Proof.} Suppose the contrary. Then there are constants $h_j$, not all \\ vanishing, such that
\begin{equation}\label{eq:8}
\sum_{j=1}^Jh_ju_{jN}(s)=0, \qquad \forall s\in S.
\end{equation}
Define $w:=\sum_{j=1}^Jh_ju_{j}(x)$. Then $w$ solves equation (\ref{eq:2}) and $w=w_N=0$ on $S$. By the uniqueness of the solution to
the Cauchy problem for equation  (\ref{eq:2}) it follows that $w=0$ in $D$.   Since
the set  $\{u_j(x)\}_{j=1}^J$ is linearly independent in $L^2(D)$, one gets $h_j=0$, $1\le j\le J$.
Lemma 2 is proved. \hfill$\Box$

{\bf Lemma 3.} {\em If  the set  $\{u_{jN}(s)\}_{j=1}^J$ is linearly independent in $L^2(S)$ and $S$ is analytic,
then  $det \left(u_{jN}(s_m)\right)\neq 0$ for almost all $s_m\in S$, $1\le m \le J. $}

{\em Proof.}  In  Lemma 3 the variables $s_m$, $1\le m \le J$, are independent. We prove that for each $s_m$
the determinant  $det \left(u_{jN}(s_m)\right)\neq 0$ for almost all $s_m\in S$, $1\le m \le J$.
Suppose the contrary, that is,  $det \left(u_{jN}(s_m)\right)= 0$ on a set $\Delta\subset S$ of positive surface measure  on $S$ for $s_m\in \Delta$, $1\le m \le J$.  Then, by analyticity of $S$ and by the resulting analyticity
of each of the elements $u_{jN}(s_m)$ for $s_m\in S$, one gets:
\begin{equation}\label{eq:9}
 det \left(u_{jN}(s_m)\right)= 0, \qquad \forall s_m\in S, \quad 1\le m \le J.
\end{equation}
Denote $s_1=s$.  If  (\ref{eq:9}) holds, then, expanding the determinant  (\ref{eq:9}) over
the elements of the first row and denoting by $A_j$ the determinant corresponding in this expansion to the element $u_{jN}(s)$, one gets:
\begin{equation}\label{eq:9'}
\sum_{j=1}^J (-1)^{j+1} u_{jN}(s) A_j=0,  \qquad \forall s\in S,
\end{equation}
 where the determinants $A_j$ do not depend on $s$. Since the set
$\{u_{jN}(s)\}_{j=1}^J$ is linearly independent in $L^2(S)$, it follows that $A_j=0$, $1\le j \le J$.
This implies that the column $\{u_{jN}(s)\}|_{j=1}^J$ is a linear combination of other columns.
Therefore, the system
$\{u_{jN}(s)\}_{j=1}^J$ is linearly dependent in $L^2(S)$, contrary to the result of Lemma 2.
This contradiction proves Lemma 3. \hfill$\Box$

\hspace{2mm}

Theorem 5 and, therefore, Theorem 3 are proved.  \hfill$\Box$

\hspace{2mm}

Let us now prove formula  (\ref{eq:4}).

\hspace{2mm}

 {\bf Theorem 6.} {\em Formula (\ref{eq:4}) holds.}

 {\em Proof.} Let $u$ solve problem (\ref{eq:2}) and let $U$ be an arbitrary solution to the
 equation $(\nabla^2+k^2)U=0$ in the ball $B_a:=\{x: |x|\le a\}$, $D\subset B_a$. Multiply
  (\ref{eq:2}) by $U$,  integrate over $D$ and use the boundary conditions (\ref{eq:2}) to get $\int_D U(x)dx=0.$
  If $g$ is an arbitrary rotation and  $(\nabla^2+k^2)U=0$, then $U(gx)$ also solves this equation. Thus,
  \begin{equation}\label{eq:12}
 \int_D U(gx)dx=0, \qquad \forall g.
\end{equation}
Let $g$ be a rotation about unit vector $\alpha$ by an angle $\phi$. Differentiate (\ref{eq:12})
with respect to $\phi$ and then let $\phi=0$. The result is:
 \begin{equation}\label{eq:13}
 \int_D \nabla U(x)\cdot [\alpha, x]dx=0, \qquad \forall \alpha\in S^2.
\end{equation}
Using the divergence theorem and arbitrariness of $\alpha\in S^2$ one obtains from (\ref{eq:13}) the relation:
 \begin{equation}\label{eq:14}
 \int_S  U(s) [s, N]ds=0.
\end{equation}
The set $\{U(s)\}$, where   $(\nabla^2+k^2)U=0$ in the ball $B_a:=\{x: |x|\le a\}$, is the orthogonal
complement in $L^2(S)$ of the linear span of the functions $u_{jN}(s)$, where $u_j(x)$, $1\le j \le J$, is a
linearly independent set of solutions to problem  (\ref{eq:5}).  To check this, denote $f:=U(s)$ and let $F\in H^2(D)$
be any function such that $F|_S=f$. Let $U=F+v$. Then $(\nabla^2+k^2)v=-(\nabla^2+k^2)F$ in $D$, $v|_S=0$.
This boundary problem for $v$ is solvable if and only if the orthogonality conditions
 \begin{equation}\label{eq:15}
\int_D  (\nabla^2+k^2)F(x) u_j(x)dx=0, \qquad 1\le j\le J,
\end{equation}
hold. Integrating by parts and taking into account that
$$u_j|_S=0 \quad and  \quad F|_S=f,$$ one gets
\begin{equation}\label{eq:16}
\int_S f(s)u_{jN}(s)ds=0, \qquad 1\le j\le J.
\end{equation}
Thus, any function orthogonal in $L^2(S)$ to $f$, that is, orthogonal  to the restriction of $U$ on $S$, is a linear
combination of the functions $u_{jN}(s)$, where the set $\{u_j(x)\}$ is a complete linearly independent set of solutions to (\ref{eq:5}).
Theorem 6 is proved. \hfill$\Box$

\hspace{3mm}

Brief historical comment: the result in \cite{P} was not correct, a counterexample was given in \cite{C}.
A bibliography on the Pompeiu problem can be found in  \cite{R363}-- \cite{R470} and in  \cite{Z}.

The method used in Theorem 6 was used in other symmetry problems, see \cite{R512}-\cite{R626}.
In \cite{R470}, p. 406,  a description is given of all non-zero solutions to equation (\ref{eq:1}).

\end{document}